\documentclass[12pt,reqno]{amsart}

\usepackage{amsmath}
\usepackage{amssymb}
\usepackage{amsthm}
\usepackage[english]{babel}
\usepackage[english]{xcolor}

\newtheorem{theorem}{Theorem}

\newcommand{\beqa}{\begin{eqnarray}}
\newcommand{\beqan}{\begin{eqnarray*}}
\newcommand{\eeqa}{\end{eqnarray}}
\newcommand{\eeqan}{\end{eqnarray*}}
\def\beq#1\eeq{\begin{equation}#1\end{equation}}

\def\P{\mathbf P }

 \def\na{\,\, {\raise.4pt\hbox{$\shortmid$}}{\hskip-2.0pt\to}\, \, }

\def\={\overset{ \text{\rm def} }=}

\def\ffrac{\frac}
\def\4{\kern1pt}

\def\bgl#1{\bigl#1\4}
\def\BR{\bf R}

\newcommand{\tc}{}
\newcommand{\gr}{}
\begin{document}

\title[\tc{\gr{Convergence} on convex polyhedra}]
{\tc{Convergence to infinite dimensional Compound Poisson Distributions on convex polyhedra}}

\author[F.~G\"otze]{Friedrich G\"otze}
\author[A.Yu. Zaitsev]{Andrei Yu. Zaitsev}

\email{goetze@math.uni-bielefeld.de}
\address{Fakult\"at f\"ur Mathematik,\newline\indent
Universit\"at Bielefeld, Postfach 100131,\newline\indent D-33501 Bielefeld,
Germany\bigskip}
\email{zaitsev@pdmi.ras.ru}
\address{St.~Petersburg Department of Steklov Mathematical Institute
\newline\indent
Fontanka 27, St.~Petersburg 191023, Russia\newline\indent
and St.Petersburg State University, 7/9 Universitetskaya nab., St. Petersburg,
199034 Russia}

\begin{abstract}
The aim of the present work is to provide a supplement to the authors' paper \cite{GZ18}. It is shown that our results
on the approximation of distributions of sums of independent summands by
the accompanying compound Poisson laws and the estimates of the proximity of sequential convolutions of multidimensional distributions on convex polyhedra may be
almost automatically transferred
to the  infinite-dimensional case.

\end{abstract}



\thanks{The authors were supported by the SFB 1283 and by RFBR-DFG grant 20-51-12004.
The second author was supported by grant RFBR 19-01-00356.}

\maketitle

The aim of the present work is to provide a supplement to the authors' paper \cite{GZ18}. It is shown that our results
on the approximation of distributions of sums of independent summands by
the accompanying compound Poisson laws and the estimates of the proximity of sequential convolutions of multidimensional distributions on convex polyhedra may be
almost automatically transferred
to the  infinite-dimensional case.
\tc{This considerably extends the scope of the results.}

Let us first introduce some notation. Let $1\le d\le\infty$ and  $\mathfrak F_d$ denote the set
of probability distributions defined on the Borel $\sigma$-field of
subsets of the Euclidean space~${\BR}^d$ and let \,${\mathcal L}(\xi)\in\mathfrak F_d$ \,be the
distribution of a $d$-dimensional random vector~$\xi$. Let
 $\mathfrak F_d^s\subset \mathfrak F_d$ be the set of symmetric distributions.
For $F\in\mathfrak F_d$, we
denote the corresponding characteristic functions by $\widehat F(t)$, $t\in \mathbf R^d$, and distribution functions
by~$F(x)=F\{(-\infty,x_1]\times\cdots\times(-\infty,x_d]\}$, $x=(x_1,\ldots,x_d)\in{\BR}^d$.
The uniform Kolmogorov distance is defined as
$$ \rho (G,H)=\sup_{x\in \mathbf R^d} \;\bgl|G(x)-H(x)\bigr|,\quad G,H\in\mathfrak F_d. $$ By the symbols~$c$ and~$c(\,\cdot\,)$ we denote (
\tc{actually various different}) positive absolute constants and quantities depending only on the arguments in brackets.
For $0\le\alpha\le2$,
we denote
$$
\mathfrak F_d^{(\alpha)}=\Big\{F\in \mathfrak F_d^s:\widehat F(t)\ge-1+\alpha,
\ \hbox{for all }t\in \mathbf R^d \Big\},\quad \mathfrak F_d^{+}=\mathfrak F_d^{(1)}.
$$Notice that our notations 
\tc{include  the case}  $d=\infty$, that is, for distributions in the Hilbert space $\mathbf R^\infty=\mathbf{H}$.

Products and powers of measures are understood in the
convolution sense: \;${GH=G*H}$, \,$H^m=H^{m*}$, \,$H^0=E=E_0$,
where $E_x$ \,is the distribution concentrated at a point $x\in
{\BR}^d $. A natural approximating infinitely divisible distribution
for $\prod_{i=1}^nF_i$ is the accompanying compound Poisson
distribution $\prod_{i=1}^n e(F_i)$, where $$
e(H)=e^{-1}\sum_{k=0}^\infty
 \ffrac{H^k}{k!},\quad H\in\mathfrak F_d,
$$ and, more generally,
\begin{equation}  e(\alpha\4 H)=e^{-\alpha}\sum_{k=0}^\infty
 \ffrac{\alpha^k\4H^k}{k!}, \quad \alpha>0.  \label{00}\end{equation}
It is well-known that the distribution $e(\alpha\4 H)$ is infinitely divisible and $e(n\4 H)=(e(H))^n$, $n \in \mathbf N$.

 Arak \cite{A80} showed that, if $F\in\mathfrak F_1^{+}$ is a symmetric one-dimensional distribution with
characteristic function which is \gr{non-negative} for all $t \in \mathbf R$, then
\begin{equation}
\rho(F^n,e(n\4F)) \le c\,n^{- 1 },
\label{a97}\end{equation}
He introduced and used the so-called method of triangular functions (see \cite[Chapter~3, Sections~2--4]{2}).

 Zaitsev \cite{z83} \tc{applied the methods 
 due to Arak} while proving inequality \eqref{a97} (see also \cite[Chapter~5, Sections~2, 5--7]{2}). 
 \tc{Later on he  modified  these methods}, adapting them to the multidimensional case (see \cite{z87}--\cite{z92}). In particular,
in \cite{z892}, a multidimensional analogue of inequality \eqref{a97} was obtained.

Using the method of triangular functions and its generalizations, several bounds of the type \begin{equation}\rho(G,H)\le c(d)\, \varepsilon\label{b97}\end{equation} were obtained, where $0<\varepsilon<1 $ is small, $G, H\in \mathfrak F_d$.

Inequality~\eqref{b97} is equivalent to the validity of the inequality
\begin{equation} \bigl|G\{P\}-H\{P\}\bigr|\le c(d)\,\varepsilon \label{g97}\end{equation}
for all sets $P$ of the form
\begin{equation} P=\big\{x\in\mathbf R^d:\langle x,e_j\rangle\le b_j, \ j=1,\ldots, d\big\}, \label{g9}\end{equation}
where $e_j\in \mathbf R^d$ are the vectors of the standard Euclidean basis, $-\infty< b_j\le \infty$, $j=1,\ldots, d$.

For $m\in\mathbf N$ we denote by $\mathcal P_m$ the collection of convex polyhedra $P \subset \mathbf R^d$
representable in the form
$$
P=\big\{x\in\mathbf R^d:\langle x,t_j\rangle\le b_j, \ j=1,\ldots, m\big\},
$$
where $t_j\in \mathbf R^d$, $-\infty\emph{}< b_j\le \infty$, $j=1,\ldots, m$.
Define also
$$ \rho_m (G,H)=\sup_{P\in\mathcal P_m} \bigl|G\{P\}-H\{P\}\bigr|. $$

Our results show that $\rho_m$ is a natural multivariate version of the Kolmogorov distance. A similar definition of the multidimensional L\'evy distance is given in \cite{GZ20}.

In the finite-dimensional case $d<\infty$, the assertions of our Theorems~\ref{Th0}--\ref{Th888} were proved in the authors' paper \cite{GZ18}, see also \cite{z80}--\cite{z92}.
\tc{We would like to emphasize  the following results remain valid for
$d=\infty$ as well.}

\begin{theorem}\label{Th0} Let $d=\infty$ and  $F\in\mathfrak F_d^{(\alpha)}$,
$0\le\alpha\le2$,
$m,n\in\mathbf N$. Then
\begin{equation} \max\big\{ \rho_m(F^n,e(nF)),  \rho_m(F^n,F^{n+1})\big\} \le c(m)\,\Big(n^{-1} +
\exp(-n\alpha+ c\,m\,\log^3 n)\Big).  \label{99}\end{equation}
\end{theorem}

It is clear that Theorem~\ref{Th0} implies inequality \eqref{a97}.

\begin{theorem}\label{Th1} Assume that $d=\infty$ and the
 distributions $G_i\in\mathfrak F_d$ are represented as
 \begin{equation}
G_i=(1-p_i)\,E+p_i\,V_i,\quad i=1,\dots,n, \label{11h}\end{equation}
where $V_i\in\mathfrak F_d$ are arbitrary distributions, $0\le p_i\le p=\max\limits_{1\le i\le n}p_i$,
$$m\in\mathbf N,\quad
 G=\prod_{i=1}^nG_i,\quad D=\prod_{i=1}^ne(G_i).$$
Then
\begin{equation} \rho_m(G,D)  \le c(m)\, p.  \label{r699}\end{equation}
\end{theorem}

The situation considered in Theorem~\ref{Th1} can be interpreted as a comparison of 
\tc{a} sample containing
\tc{independent observations of 
 rare events} (like disasters, epidemics, accidents, bankruptcies etc.) with the Poisson point process
which is obtained after a Poissonization of the initial sample (see~\cite{GZ17},~\cite{14}).
Let $Y_1, Y_2,\dots, Y_n$ be independent non-identically distributed observations taking values in  a measurable space $(\mathfrak X,\mathcal S)$. Simultaneously we observe \tc{some 
rare events}.
Assume that the set $\mathfrak X$ is represented as $\mathfrak X= \mathfrak X_1\cup\mathfrak X_2$, with \,$\mathfrak X_1,\,\mathfrak X_2\in \mathcal S$, $\mathfrak X_1\cap\mathfrak X_2=\varnothing$. We say that the $i$-th  rare event occurs if $Y_i\in\mathfrak X_2$.
Respectively, it does not occur if $Y_i\in\mathfrak X_1$.

We are interested in the cumulative loss
\begin{equation} S=f(Y_1)+
\dots+ f(Y_n), \label{639}
\end{equation}
where $f:\mathfrak X\to\mathbf  R^d$ is a loss function such that $f(x)=0$ for $x\in \mathfrak X_1\emph{}$.
Let $G_i={\mathcal L}(f(Y_i))$, $i=1,2,\dots,n$.
Then
 distributions $G_i\in\mathfrak F_d$ can be represented as mixtures \begin{equation}
G_i=(1-p_i)\,E+p_i\,V_i, \label{(11}\end{equation} where $E, V_i\in\mathfrak F_d$
are conditional distributions of vectors $f(Y_i)$ given
$Y_i\in\mathfrak X_1$ and $Y_i\in\mathfrak X_2$ respectively,
\begin{equation} 0\le p_i=\mathbf{P}\big\{Y_i\in\mathfrak X_2\big\}
=1-\mathbf{P}\big\{Y_i\in\mathfrak X_1\big\}\le1. \label{(12}\end{equation}
 \tc{We shall call events  rare} when the quantity
\begin{equation}  p=\max_{1\le i\le
n}p_i\label{(13}\end{equation}
 is small. In other words, this is   \tc{a model for uniformly rare
 events.}

 Let the random set ${\mathbf P}=\{Z_k\}$ be a realization of the Poisson point
process on the space $\mathfrak X$ with intensity measure
$\sum_{i=1}^n {\mathcal L}(Y_{i})$.
The sum $S$ defined in~\eqref{639}
has the distribution $G=\prod_{i=1}^nG_i$.
 It is easy to see that
 \begin{equation}
D=\prod_{i=1}^n e(G_i). \label{(15}
\end{equation}
 is the distribution of \begin{equation}T=\sum_{Z_k\in{\mathbf P}}f(Z_k).\label{(135}\end{equation}
 Theorem~\ref{Th1} implies that
\begin{equation}\rho_m\bigg(\mathcal{L}\Big(\sum_i f(Y_i)\Big),
\mathcal{L}\Big(\sum_{Z_k\in{\mathbf P}}f(Z_k)\Big)\bigg)\le c(m)\,p.\label{99h}\end{equation}
Thus, estimating the cumulative loss we can replace observations $Y_1, Y_2,\dots, Y_n$ by the Poisson point
process ${\mathbf P}$ if our rare events
\tc{are uniformly rare.
Then we can use the well-known  independence properties
 of disjoint sets of  the Poisson point process.}

In the rest of the paper we will study \tc{the
size of the difference between $F^{n +k}$ and $F ^n$, i.e. representing the change by addition of a finite group
of additional independent terms.}
A particular case of this problem is considered in inequality \eqref{99} of Theorem~\ref{Th0}.

 \begin{theorem}\label{Th89}
Let  $F\in\mathfrak F_d^{s}$, $d=\infty$ and
$k,m,n\in\mathbf N$. Then
\begin{equation}  \rho_m(F^n,e(nF)) \le c(m)\,n^{-1/2} ,  \label{990}\end{equation}
\begin{equation}    \rho_m(F^n,F^{n+2k})\le c(m)\,k\,n^{-1} ,  \label{b95}\end{equation}
\begin{equation}  \rho_m(F^n,F^{n+2k+1})\le  c\,m\,n^{- 1/2 }+c(m)\,k\,n^{-1} .  \label{c903}\end{equation}In particular,
\begin{equation}  \sup_{k\le\sqrt{n}} \rho_m(F^n,F^{n+k})\le  c(m)\,n^{- 1/2 } .  \label{c904}\end{equation}
\end{theorem}

 For $m\in\mathbf N$, $t_1,\ldots,t_m\in\mathbf R^d$, we denote by $\mathcal P(t_1,\ldots,t_m)$ the collection of convex polyhedra $P \subset \mathbf R^d$
representable in the form
$$
P=\big\{x\in\mathbf R^d:\langle x,t_j\rangle\le b_j, \ j=1,\ldots, m\big\},
$$
where  $-\infty\emph{}< b_j\le \infty$, $j=1,\ldots, m$.
 Clearly,
 $$
\mathcal P_m=\bigcup_{t_1,\ldots,t_m}\mathcal P(t_1,\ldots,t_m).
$$

\begin{theorem}\label{Th88} Let  $F=\mathcal{L}(\xi)\in\mathfrak F_d$, $d=\infty$, $m\in\mathbf N$, $t_1,\ldots,t_m\in\mathbf R^d$, and  all distributions of the random variables $\langle \xi,t_j\rangle$, $j=1,\ldots, m$, are either non-degenerate
or equal to~$E\in \mathfrak F_1$. Then,
for all
$P\in\mathcal P(t_1,\ldots,t_m)$,
\begin{equation} \bigl| (F^n) \{P\} - (F^{n+1}) \{P\}\big|  \le c(F, t_1,\ldots,t_m)\,n^{-1/2} .  \label{333}\end{equation}
\end{theorem}

 It is evident that if a distribution $F\in\mathfrak F_d$
is concentrated on a hyperplane which does not contain zero, then
$\rho_m(F^n,F^{n+k}) =1$ for any $m,n, k\in\mathbf N$.
Thus, we have \gr{arrived} \tc{at  the following alternative: either the left-hand side of \eqref{333} is equal to one or it decreases at least as $O(n^{-1/2})$.}

In conclusion, we formulate a result on the 
\gr{prozimity} of the distributions of sums of random number of independent identically distributed random vectors, which follows from Theorem~\ref{Th89}.
For the distance $\rho(\,\cdot\,,\,\cdot\,)$ this result is contained in \cite[Theorem 1.3]{z88}.

Let $\xi_1, \xi_2, \ldots$ be independent identically distributed random vectors with common distribution $F\in\mathfrak F_d$ and let $(\mu, \nu)\in\mathbf{Z}^2$ be a two-dimensional random vector with integer non-negative coordinates, independent of the sequence $\{\xi_j\}_{j=1}^\infty$. Denote
\begin{equation}
U=\mathcal{L}(\mu), \enskip V=\mathcal{L}(\nu), \enskip G=\mathcal{L}(\xi_1+\cdots+\xi_\mu)
, \enskip H=\mathcal{L}(\xi_1+\cdots+\xi_\nu).
  \label{007a}\end{equation}
Then it is well known that
\begin{equation}
G=\sum_{k=0}^\infty\P\{\mu=k\}\,F^k,\quad H=\sum_{k=0}^\infty\P\{\nu=k\}\,F^k.
  \label{008a}\end{equation}

\begin{theorem}\label{Th888}Let $d=\infty$. If $F\in\mathfrak F_d^s$, then
\begin{equation}
\rho_m(G,H) \le\inf{\mathbf E}\,\min\Big\{\frac{c\4m}{\sqrt{\nu+1}}+c(m)\,\frac{\left|\mu-\nu\right|}{\nu+1}, 1\Big\},
  \label{708a}\end{equation}
and if $F\in\mathfrak F_d^+$, then
\begin{equation}
\rho_m(G,H) \le\inf{\mathbf E}\,\min\Big\{c(m)\,\frac{\left|\mu-\nu\right|}{\nu+1}, 1\Big\}.
  \label{608a}\end{equation}
Here, the infimum is taken over all possible two-dimensional distributions $\mathcal{L}((\mu, \nu))\in\mathfrak F_2$ such that  $\mathcal{L}(\mu)=U$, $\mathcal{L}(\nu)=V$.
\end{theorem}

\noindent {\it Proof of Theorems\/ $\ref{Th0}$--$\ref{Th888}$.}
 At first \tc{note} that Theorems \ref{Th0}--\ref{Th888} for $d<\infty$ were proved in \cite{GZ18}.

 Fix some  polyhedron $P\in\mathcal P_m$:
\begin{align*}
    P=\big\{x\in\mathbf H: \langle x,t_j\rangle\le b_j, \ j=1,\ldots, m\big\},
\end{align*}where $t_j\in\mathbf H$, $\|t_j\|=1$, $b_j\in\mathbf R$,
$j=1,\ldots, m$. Let $\mathbf L_t\subset\mathbf H$ be the linear span of vectors
$$
\big\{t_j, j=1,\ldots, m\big\}, \quad k=\dim \mathbf L_t\le m,
$$
and let $\mathbb P_t:\mathbf H\to\mathbf L_t$ be the orthogonal projection operator on the subspace $\mathbf L_t$.
Consider the polyhedron $\overline{P}\subset\mathbf L_t$  defined as
\begin{align*}
 \overline   P=\big\{x\in\mathbf L_t: \langle x,t_j\rangle\le b_j, \ j=1,\ldots, m\big\}.
\end{align*}
It is easy to see that, for any random vector $\zeta\in\mathbf H$, we have
\begin{align*}
  \langle\mathbb   P_t\zeta,t_j\rangle=\langle \zeta,t_j\rangle,  \quad j=1,\ldots, m.
\end{align*}
Therefore,
 \begin{equation}
    \P\{\zeta\in P\}=  \P\{\mathbb P_t\zeta\in \overline{P}\}.
\label{div7}
\end{equation}
The distributions of $k$-variate vectors  $\mathbb P_t \xi$, $\mathbb P_t \eta$  satisfy the same $k$-dimensional conditions as the distributions of the random vectors $\xi,\eta\in\mathbf H$ with compared infinite-dimensional distributions. For example, if
$\mathcal{L}(\xi)\in\mathfrak F_\infty^{(\alpha)}$,
for some $\alpha$ satisfying $0\le\alpha\le2$, then $\mathcal{L}\big(\mathbb P_t\xi\big)\in\mathfrak F_k^{(\alpha)}$.
Similarly, if
$\mathcal{L}(\xi)\in\mathfrak F_\infty^{s}$, then $\mathcal{L}\big(\mathbb P_t\xi\big)\in\mathfrak F_k^{s}$ and so on. It remains to apply the $k$-dimensional versions of Theorems \ref{Th0}--\ref{Th888}.\hfill$\square$\medskip

An interesting problem is to \tc{extend} our inequalities to arbitrary convex sets~$P$. This does not follow from Theorems \ref{Th0}--\ref{Th888}, since the constants $c(m)$ depend on~$m$.

\end{document}